\documentclass{article}

\usepackage{natbib}
\usepackage{amsmath,amsfonts,amssymb}
\usepackage[arrow, matrix, curve]{xy}
\usepackage{enumerate}

\begin{document}

\title{Implicitization of rational ruled surfaces with $\mu$-bases}

\author{Marc Dohm}
\maketitle

\begin{abstract}
\noindent Chen, Sederberg, and Zheng introduced the notion of a
$\mu$-basis for a rational ruled surface in \cite{MR1807798} and
showed that its resultant is the implicit equation of the surface,
if the parametrization is generically injective. We generalize
this result to the case of an arbitrary parametrization of a rational
ruled surface. We also 
give a new proof for the corresponding theorem in the curve case
and treat the reparametrization problem for curves and ruled
surfaces. In particular, we propose a partial solution to the problem of computing a proper re\-pa\-ra\-me\-tri\-za\-tion for a rational ruled surface.
\end{abstract}

\section*{Introduction}

\noindent Implicitization is a fundamental problem in Computer
Aided Geometric Design and there are numerous applications related
to it, e.g. the computation of the intersection of two ruled
surfaces, see \cite{GonzVega}. The method of $\mu$-bases (also
known as ``moving lines" or ``moving surfaces") constitutes an
efficient solution to the implicitization problem. Introduced in
1998 by Cox, Sederberg, and Chen for rational curves in
\cite{MR1638732}, it was generalized to ruled surfaces in
\cite{MR1807798} and \cite{MR2015108}. Whereas the curve case is
very well understood and we know that the resultant of a
$\mu$-basis is the implicit equation to the power $d$, where $d$
is the degree of the rational map induced by the parametrization, this result is still to be
shown in its full generality (i.e. for arbitrary $d$) for ruled
surfaces. We fill this gap by giving a proof, which relies on a
geometric idea that reduces the ruled surface case to the curve
case. From a computional point of view, $\mu$-bases are in general
more efficient than other resultant-based methods such as the ones
introduced in \cite{MR2172855} or in \cite{MR2004036}, since they
are well adapted to the geometry of ruled surfaces.

\section{$\mu$-bases of rational planar curves}

\noindent As we will need them later on, we will start with some
known results about the $\mu$-basis of a rational parametric
planar curve $\mathcal{C}$ over an algebraically closed field
$\mathbb{K}$ of arbitrary characteristic, i.e. one given by a parametrization map
$$\begin{array}{rcl}
    \Phi_\mathcal{C} : \quad \mathbb{P}^{1} & \dashrightarrow&  \mathbb{P}^{2}
      \nonumber \\
     (s:\bar{s}) &\mapsto& (f_0(s,\bar{s}):f_1(s,\bar{s}):f_2(s,\bar{s}))
\end{array}$$
where each $f_{i} \in \mathbb{K}[s,\bar{s}] =: R$ is homogeneous
of degree $n>0$ and $g := gcd(f_0,f_1,f_2)$ is of degree strictly less than
$n$. The first syzygy module of $f_0,f_1,f_2$ is defined as
$$\textsl{Syz}(f_{0},f_{1},f_{2})=\{ P \in R[x,y,z] \: | \: deg(P)\leq 1, P(f_{0},f_{1},f_{2}) = 0 \} \subseteq
R[x,y,z]$$ Then we have the following well-known result.

\newtheorem{hst}{Theorem}
\begin{hst}
There exists an isomorphism of graded $R$-modules $$
\textsl{Syz}(f_{0},f_{1},f_{2}) \cong R(-\mu_1) \oplus R(-\mu_2)$$
where $\mu_i \in \mathbb{N}$, $\mu_1 \leq  \mu_2$ and
$$\mu_1 +\mu_2= n-deg(g) = deg(\Phi_\mathcal{C}) \cdot
deg(\mathcal{C})=: d$$   \label{syzisocurv}
\end{hst}

\noindent The isomorphism in the above theorem is a direct consequence
of the Hilbert-Burch Theorem (see \cite{MR1322960}[Th. 20.15]) applied to the 
exact sequence 
\begin{equation*} \label{seq1} 0 \rightarrow \textsl{Syz}(f_0,f_1,f_2)(-n)
\rightarrow R^{3}(-n) \rightarrow R \rightarrow R/I \rightarrow 0 \end{equation*}
and the degree property can easily be checked by computing the Hilbert polynomials of this sequence.\\

\noindent A basis $(p,q)$ of $\textsl{Syz}(f_{0},f_{1},f_{2})$
with minimal degrees $deg(p)=\mu_1$ and $deg(q)=\mu_2$ in $s$ and
$\bar{s}$ is called a $\mu$-basis of the parametrization
$\Phi_\mathcal{C}$. One interesting feature of $\mu$-bases is that
the resultant of its elements is a power of the implicit equation
of $\mathcal{C}$, as was proved in \cite[Sect. 4, Th.
1]{MR1638732}. We propose an alternative proof which relies on the
idea that we can reduce the problem to the generically injective case.
The essential tool for this reduction is the existence of a proper reparametrization, which
is a consequence of L\"{u}roth's Theorem, a proof
of which can be found for example in \cite[Section
5.4]{MR0263582}. In the following lemma we deduce a reparametrization with an additional property.

\newtheorem{reparacurv}[hst]{Lemma}
\begin{reparacurv}
There exists $\psi : \mathbb{P}^{1} \dashrightarrow
\mathbb{P}^{1}$ parametrized by two coprime homogeneous
polynomials $h_0$ and $h_1$ of degree $deg(\Phi_\mathcal{C})$ and
a parametrization $\Phi'$ of $\mathcal{C}$ defined by homogeneous
polynomials $f'_0(s,\bar{s})$,$\; f'_1(s,\bar{s})$ and
$f'_2(s,\bar{s})$ such that the following diagram commutes:

$$\begin{xy}
 \xymatrix{
    \mathbb{P}^1 \ar@{-->}[rrr]^{\Phi_\mathcal{C}} \ar@{-->}[dd]_\psi & & & \mathbb{P}^2  \\
                                            & & &               \\
    \mathbb{P}^1 \ar@{-->}[rrruu]_{\Phi'_\mathcal{C}}         & & &
    }
 \end{xy}$$ \\
It follows that $\Phi'_\mathcal{C}$ is a proper (i.e. generically
injective) parametrization of $\mathcal{C}$, in other words
$deg(\Phi'_\mathcal{C})~=~1$. Moreover, if
$gcd(f_0,f_2)=gcd(f_1,f_2)=1$, we can choose $\Phi'_\mathcal{C}$
such that $f_i =f'_i(h_0,h_1)$ for $i \in \{0,1,2\}$. \\
\label{repacu}
\end{reparacurv}

\noindent \textit{Proof:} First, we treat the case
$gcd(f_0,f_2)=gcd(f_1,f_2)=1$. Then we can dehomogenize
$\frac{f_0}{f_2}$ and $\frac{f_1}{f_2}$ by setting $\bar{s}=1$
without changing the degree as rational functions and decompose
them by means of L\"{u}roth's Theorem \cite[Section
5.4]{MR0263582}) in the following way $$\begin{array}{rcl}
\frac{f_0}{f_2} = \frac{f'_0}{f'_2} \circ \frac{h_0}{h_1} & \qquad
\qquad \qquad & \frac{f_1}{f_2} = \frac{f'_1}{\tilde{f}'_2} \circ
\frac{h_0}{h_1}
\end{array}$$ with $gcd(h_0,h_1)=gcd(f'_0,f'_2)=gcd(f'_1,\tilde{f}'_2)=1$ and
$deg(h_0)=deg(h_1)=deg(\Phi_\mathcal{C})$ after having
rehomogenized them with respect to $\bar{s}$. By multiplying the
fractions with a suitable power of $h_1$ we can consider the
$f'_i$ as bivariate homogeneous polynomials
$$\begin{array}{rcl} \frac{f_0}{f_2} =
\frac{f'_0(h_0,h_1)}{f'_2(h_0,h_1)} & \qquad \qquad \qquad &
\frac{f_1}{f_2} = \frac{f'_1(h_0,h_1)}{\tilde{f}'_2(h_0,h_1)}
 \nonumber
\end{array}$$ 
Then the numerators and denominators are all coprime,
which for the right hand sides follows from \cite[Prop.
6]{ZippelRFD} and we deduce the term-by-term equalities $f_i
=f'_i(h_0,h_1)$ for $i \in \{0,1,2\}$. \\
\noindent In the general case, we divide the polynomials of the
parametrization by their greatest common divisor and perform a
generic coordinate change in order to pass to another
parametrization of $\mathcal{C}$ which fulfills
$gcd(f_0,f_2)=gcd(f_1,f_2)=1$ and whose polynomial decomposition
completes the commutative diagram
of rational maps. $\square$ \\

\noindent Now we are ready to proceed to the main theorem of this
section, for which we give a new proof that establishes a link
between the $\mu$-basis of $\Phi_\mathcal{C}$ and a $\mu$-basis of
a proper reparametrization of the curve.

\newtheorem{impleqcurv}[hst]{Theorem}
\begin{impleqcurv}
Let $(p,q)$ be a $\mu$-basis of the parametrization
$\Phi_\mathcal{C} : \mathbb{P}^{1} \dashrightarrow
\mathbb{P}^{2}$. Then
$$Res(p,q)=F_\mathcal{C}^{deg(\Phi_\mathcal{C})}$$ where $F_\mathcal{C}$ is an implicit
equation of the curve $\mathcal{C}$ defined by $\Phi_\mathcal{C}$
and $Res(p,q) \in \mathbb{K}[x,y,z]$ is the homogeneous resultant
with respect to the indeterminates $s$ and $\bar{s}$.
\label{implcurv} \\
\end{impleqcurv}

\noindent \textit{Proof:} First of all, we may assume that
$gcd(f_0,f_2)=gcd(f_1,f_2)=1$ (if necessary, we divide by
$gcd(f_0,f_1,f_2)$ and perform a generic coordinate change, both
of which do not affect the result). So by Lemma
\ref{repacu} there exist $f'_0,f'_1,f'_2 \in R$ and homogeneous,
coprime $h_0, h_1 \in R$ of degree $deg(\Phi_\mathcal{C})$, such
that $$\begin{array}{rcl}
f_0 &=&f'_0(h_0,h_1)  \\
f_1 &=&f'_1(h_0,h_1)\\
f_2 &=&f'_2(h_0,h_1)
\end{array}$$
\noindent Let $(p',q')$ be a $\mu$-basis of the proper
reparametrization $\Phi'_\mathcal{C}$ of $\mathcal{C}$ defined by
the $f'_i$. Then $p'(h_0,h_1)$ and $q'(h_0,h_1)$ are linearly
independent syzygies (i.e. we substitute $h_0$ for $s$ and $h_1$
for $\bar{s}$). It is easy to see that they form a $\mu$-basis by
verifying the degree property and if $\mu_1 < \mu_2$, they are
related to our original $\mu$-basis $(p,q)$ by
$$\begin{array}{rcl}
 p'(h_0,h_1) &=& \lambda p  \\
 q'(h_0,h_1) &=& ap + q
 \end{array}$$
 for some constant $\lambda \neq 0$ and a homogeneous $a \in R$ of
 degree $deg(q)-deg(p)$. (If $\mu_1=\mu_2$, we have $p' \circ h=\alpha_1 p +\alpha_2
 q$ and $q' \circ h = \beta_1 p + \beta_2 q$ for some constants
 $\alpha_i$ and $\beta_i$ (see \cite[Th. 2]{ChenWangCurves}), which leads to computations that are analogous to the ones
 that follow).

 \noindent Now we can apply elementary properties of resultants to calculate
 \begin{align}
 Res(p,q) &= \lambda^{-\mu_2} \cdot Res(\lambda p,ap + q)   \nonumber \\
          &= \lambda^{-\mu_2} \cdot Res(p'(h_0,h_1),q'(h_0,h_1))  &        \nonumber \\
          &= \lambda^{-\mu_2} \cdot Res(h_0,h_1)^{deg(p')deg(q')} \cdot
          Res(p',q')^{deg(h_0)}    & \nonumber \\
          &= c \cdot Res(p',q')^{deg(\Phi_\mathcal{C})} & (c \in \mathbb{K}^\ast)      \label{resrechnung}
 \end{align}
where $c=\lambda^{-\mu_2} \cdot Res(h_0,h_1)$ is a constant
(since the $h_i$ do not
depend on $x,y,z$) and non-zero (because $gcd(h_0,h_1)=1$). The third identity is a well-known base change
formula for resultants, which is proved in
\cite[5.12]{JouanolouRes}, and in the last identity we used $deg(h_0)=deg(\Phi_\mathcal{C})$.

\noindent So by \eqref{resrechnung} we have reduced the theorem to the
special case where the parametrization has degree 1, and it
remains to show: $$\begin{array}{rcl}
a) & \quad & Res(p',q') \neq 0  \nonumber \\
b) & \quad & F_\mathcal{C} \: | \: Res(p',q') \nonumber \\
c) & \quad & deg_{x,y,z}(Res(p',q')) \leq deg(\mathcal{C})
\nonumber
\end{array}$$

\begin{enumerate}[a)]
\item Suppose $p= G \cdot H$ were reducible into
non-constant $G,H \in R[x,y,z]$, then one of the two, say $G$,
would be independent of $x,y,z$, because $p$ is linear in those
variables and $H$ would define a syzygy with lower degree than $p$
which contradicts the definition of a $\mu$-basis. So $p$ is
irreducible in $R[x,y,z]$ and $Res(p,q)=0$ would mean that $q=r \cdot p$ with $r \in R$,
which is impossible, for $p$ and $q$ are linearly independent over $R$. Hence $Res(p,q)\neq 0$ and by \eqref{resrechnung}
also $Res(p',q') \neq 0$.
\item By construction $p$ and $q$ vanish for all
points in $Im(\Phi_\mathcal{C})$. So for any $X=(x_1:x_2:x_3) \in
Im(\Phi_\mathcal{C})$ we have that $p(X)=q(X)=0$ which rests true
after setting $\bar{s}=1$, so the two univariate polynomials have
a common zero and therefore $Res(p(X),q(X))=(Res(p,q))(X)=0$.
Again, by \eqref{resrechnung} we have $(Res(p',q'))(X)=0$ as well and
it follows that the implicit equation $F_\mathcal{C}$ divides $Res(p',q')$.
\item All the coefficients of $p$ and $q$ are of
degree $\leq 1$ in $x,y,z$, so we can give an upper bound for the
degree of the resultant in $x,y,z$: $$deg_{x,y,z}(Res(p,q)) \leq deg(p)+deg(q)=d=deg(\Phi_\mathcal{C})deg(\mathcal{C})$$ Once again
we look at \eqref{resrechnung} to deduce that $deg_{x,y,z}(Res(p',q')) \leq deg(\mathcal{C})$ which
concludes the proof. \nolinebreak $\square$
\end{enumerate}

\section{Implicitization of rational ruled surfaces with $\mu$-bases}

\noindent Chen, Sederberg, and Zheng introduced the notion of a
$\mu$-basis for rational ruled surfaces in \cite{MR1807798}, and
it was further developed in \cite{MR2015108}. However, they worked
with the restrictive assumption that the parametrization is
generically injective. In this section, we will give a proof for
the ruled surface version of Theorem \ref{implcurv} in its general
form and explain to what extent the ruled surface case can be
reduced to the curve case.

\noindent In this paper, a rational ruled surface $\mathcal{S}$ is meant to be
a surface given by a rational map $$\begin{array}{rcl}
    \Phi_\mathcal{S} : \quad \qquad \mathbb{P}^{1} \times \mathbb{P}^{1}& \dashrightarrow&  \mathbb{P}^{3}
      \nonumber \\
     ((s:\bar{s}),(t:\bar{t})) &\mapsto& (f_0(s,\bar{s},t,\bar{t}):\ldots
     :f_3(s,\bar{s},t,\bar{t})) \label{eqsurfdef}
\end{array}$$
where  the $f_i \in \mathbb{K}[s,\bar{s},t,\bar{t}]$ are
bihomogeneous of degree $(n,1)$, by which we mean that they are
homogeneous of degree $n+1$ and that $deg_{s,\bar{s}}(f_i)=n$ and
$deg_{t,\bar{t}}(f_i)=1$ for all $i=0,\ldots,3$. We assume that
$gcd(f_0,\ldots,f_3)=1$ and that we can rewrite
\begin{equation} f_i = \bar{t}\bar{s}^{n_1-n_0}f_{i0}+t f_{i1}
\label{eqhompolys}
\end{equation}
where $f_{i0},f_{i1} \in \mathbb{K}[s,\bar{s}]$ and
$n_0:=max(deg_s(f_{i0}))$ and $n_1:=max(deg_s(f_{i1}))$, and where
we have assumed that $n_1 \geq n_0$ (otherwise we may
reparametrize \eqref{eqsurfdef} by exchanging $t$ and $\bar{t}$)
and $n_1=n$ (otherwise, we may divide the $f_i$ by a suitable
power of $\bar{s}$). Finally, we need to make the assumption that
$(f_{00},\ldots,f_{30})$ and $(f_{01},\ldots,f_{31})$ are
$R$-linearly independent to exclude the degenerate case where
$\Phi_\mathcal{S}$ does not parametrize a surface.\\
Let us fix some notation first: The $R$-module of syzygies on
$f_{0},\ldots,f_3$ depending only on $s$ and $\bar{s}$ is defined
as
$$\textsl{Syz}_R(f_{0},\ldots,f_{3})=\{ P \in R[x,y,z,w] \: | \: deg(P) = 1,
P(f_{0},f_{1},f_{2},f_{3}) = 0 \}$$ Then the structure of this
module is well known; see
\cite{MR1807798} for a proof of the following

\newtheorem{syzsurfacetheo}[hst]{Theorem}
\begin{syzsurfacetheo}
There exists an isomorphism of graded $R$-modules $$
\textsl{Syz}_{R}(f_{0},\ldots,f_{3}) \cong R(-\mu_1) \oplus
R(-\mu_2)$$ where $\mu_i \in \mathbb{N}$, $\mu_1 \leq \mu_{2}$ and
$\mu_1+\mu_2 = deg(\Phi_\mathcal{S}) \cdot deg(\mathcal{S}).$
\label{syzisosurf} \\
\end{syzsurfacetheo}

\noindent A basis $(q_1,q_2)$ of ${Syz}_R(f_{0},f_{1},f_{2},f_3)$
where $q_1$ and $q_2$ are homogeneous of minimal degrees
$deg(q_1)=\mu_1$ and $deg(q_2)=\mu_2$ in $s$ and $\bar{s}$ is
called a $\mu$-basis of the parametrization $\Phi_\mathcal{S}$. As
we can see, the syzygy module of the surface $\mathcal{S}$
resembles the one of a curve, which leads to the following question: is there 
a curve with the same syzygy module which can be defined by
means of the surface parametrization? The answer to this question
is positive and according to an idea due to \cite{BEG05}, we define the curve $\mathcal{C}$ associated to
$\mathcal{S}$ by $$\begin{array}{rcl}
    \Phi_\mathcal{C} : \quad \mathbb{P}^{1} & \dashrightarrow&  \mathbb{P}^{2}
      \nonumber \\
     (s:\bar{s}) &\mapsto& (p_{03}(s,\bar{s}):p_{13}(s,\bar{s}):p_{23}(s,\bar{s})) \nonumber
\end{array}$$
where $p_{ij} := f_{i0}f_{j1}-f_{i1}f_{j0} \in R$ are the
Pl\"{u}cker coordinates, which are homogeneous of degree
$n_1+n_0$. Let us denote $g :=gcd(p_{03},p_{13},p_{23})$.\\
\noindent The geometric idea behind this definition is that for
almost all parameter values $(s:\bar{s}) \in \mathbb{P}^1$ the
image of the map
$$\begin{array}{rcl}
    \Phi_\mathcal{S}((s:\bar{s}),-) : \quad \qquad \mathbb{P}^{1} & \dashrightarrow&  \mathbb{P}^{3}
      \nonumber \\
     (t:\bar{t}) &\mapsto& (f_0(s,\bar{s},t,\bar{t}):\ldots
     :f_3(s,\bar{s},t,\bar{t})) 
\end{array}$$ is a line $L_{(s:\bar{s})}$ in
$\mathbb{P}^3$, hence the surface $\mathcal{S}$ can be viewed as
the closure of the union of these lines. The curve defined by all
the Pl\"{u}cker coordinates $$\begin{array}{rcl}
    \Psi : \quad  \mathbb{P}^{1} & \dashrightarrow&  \mathbb{P}^{5}
      \nonumber \\
     (s:\bar{s}) &\mapsto& (p_{ij})_{i,j \in \{0,\ldots,3\},\; i<j}  \nonumber
\end{array}$$
is contained in a quadric parametrizing the lines in
$\mathbb{P}^3$, more precisely there is a one-to-one
correspondance between the points $\Psi((s:\bar{s}))$ on the
Pl\"{u}cker curve and the lines $L_{(s:\bar{s})}$ on the ruled
surface $\mathcal{S}$, which will allow us to carry over the
results about curves to the ruled surface case. However, it is
more convenient to work with the curve $\Phi_\mathcal{C}$, which
is a projection of $\Psi$ to $\mathbb{P}^{2}$. As we will see, we
need to make sure that this projection does not add any base
points, which is the statement of the following lemma.

\newtheorem{plueckerprojection}[hst]{Lemma}
\begin{plueckerprojection}
If $gcd(f_{30},f_{31})=1$ then
$$gcd(p_{03},p_{13},p_{23})=gcd(p_{03},p_{13},p_{23},p_{01},p_{02},p_{12})$$
\end{plueckerprojection}

\noindent \textit{Proof:} Let us suppose
$q=gcd(p_{03},p_{13},p_{23}) \neq 1$; the case $q=1$ is trivial.
We need to show that $q$ divides the other Pl\"{u}cker coordinates
as well. Euclidean division of the $f_{ij}$ by $q$ yields
$$f_{ij}=q \cdot \tilde{f}_{ij} + a_{ij}$$ We have the congruences
$$ p_{ij}\equiv f_{i0}f_{j1}-f_{i1}f_{j0}\equiv
a_{i0}a_{j1}-a_{i1}a_{j0} \qquad (mod \: q)$$ The other cases
being analogous, we only show $p_{12} \equiv 0 \quad (mod \; q)$,
i.e. that $a_{10}a_{21}-a_{11}a_{20}$ is divisible by $q$. Since $p_{13}$ and $p_{23}$ are
divisible by $q$, we can
write $a_{10}a_{31}-a_{11}a_{30}=qr_1$ and
$a_{20}a_{31}-a_{21}a_{30}=qr_2$, or equivalently $a_{21}a_{30}=a_{20}a_{31}-qr_2$ and $a_{11}a_{30}=a_{10}a_{31}-qr_1$. As $gcd(f_{30},f_{31})=1$ it follows that not both $f_{30}$ and $f_{31}$ are divisible by $q$, so we may assume that one of the rests of the Euclidean division, say $a_{30}$, is non-zero.
We have $$a_{30}(a_{10}a_{21}-a_{11}a_{20})=a_{10}(a_{20}a_{31}-qr_2)-a_{20}(a_{10}a_{31}-qr_1)=q
\cdot (r_1-r_2)$$ and as $a_{30}$ is non-zero and prime to $q$, we
conclude that $a_{10}a_{21}-a_{11}a_{20}$ is divisible by $q$.
$\square$ \\

\noindent Later on, we will see in another context why the
condition $gcd(f_{30},f_{31})=1$ is necessary. We should note that
it is non-restrictive, since it can always be achieved by a
generic coordinate change. Next, we state a useful degree formula,
which we will use to study the relationship between a ruled
surface and its associated curve in more detail.

\newtheorem{degformula}[hst]{Proposition}
\begin{degformula}[Degree Formula]
With the same notation and hypotheses as before the equality
$$deg(\mathcal{S})deg(\Phi_\mathcal{S})=n_1+n_0-deg(g)$$
holds. \\
\end{degformula}

\noindent \textit{Proof:} This formula is an adaptation of the
general result
$$deg(\mathcal{S})deg(\Phi_\mathcal{S})=2n - \sum_{p \in V(f_0,\ldots,f_3)} m_p$$ (see
\cite[Prop. 4.4]{FultonInter} for a proof, $m_p$ is the
multiplicity of $p$). Our formula follows by counting the base
points $\sum_{p \in V(I)} m_p=deg(g)+(n_1-n_0)$, where $n_1-n_0$
is the trivial multiplicity of the base point
$(\infty,0):=((1:0),(0:1))$ and where the other base points
(including additional multiplicities of $(\infty,0)$) can be
identified with the roots of $g$ by elementary
calculations. $\square$ \\

\noindent Note that for characteristic zero $deg(\Phi_\mathcal{S})$ - and thus also
$deg(\mathcal{S})$ - can be computed by means of 
gcd and resultant computations, see \cite{MR2076380}. 

Next, we proceed to relate
$\textsl{Syz}_{R}(f_{0},\ldots,f_{3})$ to the syzygy module of the
associated curve, given as
$$\textsl{Syz}(p_{03},p_{13},p_{23})=\{ P \in R[x,y,z] \: | \:
deg(P)= 1, P(p_{03},p_{13},p_{23}) = 0 \} $$
\newtheorem{syzsurfcurv}[hst]{Proposition}
\begin{syzsurfcurv}
If $gcd(f_{30},f_{31})=1$, then there exists a canonical
isomorphism of graded $R$-modules
$$ \textsl{Syz}_{R}(f_{0},\ldots,f_{3}) \cong
\textsl{Syz}(p_{03},p_{13},p_{23})$$  and $deg(\Phi_\mathcal{S})
\cdot deg(\mathcal{S})=deg(\Phi_\mathcal{C}) \cdot
deg(\mathcal{C}).$ \label{isofsyz}\\
\end{syzsurfcurv}

\noindent \textit{Proof:} As a direct consequence of Theorem
\ref{syzisocurv} and the degree formula, we obtain
$$deg(\Phi_\mathcal{C}) \cdot
deg(\mathcal{C})=n_1+n_0-deg(g)=deg(\Phi_\mathcal{S}) \cdot
deg(\mathcal{S})$$ and it remains to construct an isomorphism of
degree zero between the syzygy modules. Let $h_0x+h_1y+h_2z+h_3w
\in \textsl{Syz}_{R}(f_{0},\ldots,f_{3})$. As it does not depend
on $t$ and $\bar{t}$, we can deduce from \eqref{eqhompolys} that
$$\begin{array}{rcl}
h_0f_{00}+h_1f_{10}+h_2f_{20}+h_3f_{30} &=&0 \nonumber\\
h_0f_{01}+h_1f_{11}+h_2f_{21}+h_3f_{31} &=&0 \nonumber
\end{array}$$
By multiplying the first equation by $f_{31}$ and the second one
by $f_{30}$ and by substracting the second from the first we get
\begin{equation}
h_0p_{03}+h_{1}p_{13}+h_{2}p_{23}=0 \label{eqpluecksyz}
\end{equation}
which is a syzygy on the $p_{i3}$. Hence, by setting $w=0$ we
obtain a well-defined morphism $$\begin{array}{rcl}
    \varphi : \qquad \quad \textsl{Syz}_{R}(f_{0},\ldots,f_{3}) & \rightarrow&  \textsl{Syz}(p_{03},p_{13},p_{23})
      \nonumber \\
     h_0x+h_1y+h_2z+h_3w &\mapsto& h_0x+h_1y+h_2z  \nonumber
\end{array}$$
which has obviously degree zero. Now $\varphi$ is injective, because if $h_0=h_1=h_{2}=0$ for a
syzygy on the $f_i$, then $h_3=0$ as well (as $f_{30}$ and $f_{31}$ are coprime
and hence non-zero). To see why it is also
surjective, let $h_0x+h_1y+h_2z \in
\textsl{Syz}(p_{03},p_{13},p_{23})$ and by rewriting
\eqref{eqpluecksyz} we have
$$(h_0f_{00}+h_{1}f_{10}+h_{2}f_{20})f_{31}=(h_0f_{01}+h_1f_{11}+h_{2}f_{21})f_{30}$$
The assumption that $f_{30}$ and $f_{31}$ are coprime implies that
there is a polynomial $h \in K[s,\bar{s}]$ such that
\begin{equation}
 hf_{30}=h_0f_{00}+h_{1}f_{10}+h_{2}f_{20} \label{inversesyz}
\end{equation}
 and by substituting this in the above equation also
$hf_{31}=h_0f_{01}+h_1f_{11}+h_{2}f_{21}$. These two relations
show that $h_0x+h_1y+h_2z-hw \in
\textsl{Syz}_{R}(f_{0},\ldots,f_{m})$ is a preimage of $h_0x+h_1y+h_2z$, hence $\varphi$ is
surjective and the proof is complete. \nolinebreak $\square$  

\newtheorem{equaldegrees}[hst]{Corollary}
\begin{equaldegrees}
If we perform a generic coordinate change beforehand, we also have
$deg(\mathcal{S})=deg(\mathcal{C})$ and
$deg(\Phi_\mathcal{S})=deg(\Phi_\mathcal{C})$ in the situation of
the preceding Proposition \ref{isofsyz}.
\end{equaldegrees}

\noindent \textit{Proof:} As we have seen in the proof of the
proposition, the associated curve is obtained by intersecting the
surface with the plane $w=0$ and the isomorphism of the syzygy
modules is induced by the projection map. If this plane is
generic, the theorem of B\'{e}zout ensures that this intersection
preserves the degree. $\square$ \\

\noindent \noindent  An important remark is that the inverse of $\varphi$ in the proof of Proposition \ref{isofsyz}  can be described explicitly as
\begin{eqnarray}
    \varphi^{-1} : \quad  \textsl{Syz}(p_{03},p_{13},p_{23}) & \rightarrow&  \textsl{Syz}_{R}(f_{0},\ldots,f_{3})
      \label{syzisom} \\
     h_0x+h_1y+h_2z  &\mapsto&  h_0x+h_1y+h_2z-\frac{h_0f_{00}+h_{1}f_{10}+h_{2}f_{20}}{f_{30}}w  \nonumber
\end{eqnarray} by using equation \eqref{inversesyz}.It is of degree 0 and
hence preserves degrees, so it takes $\mu$-bases to $\mu$-bases. This leads to an
efficient method for the computation of the $\mu$-basis of the
surface: One computes the $\mu$-basis of the associated curve
and takes its image under $\varphi^{-1}$. See Section \ref{algoex} for
an explicit description of this algorithm.

One can regard the results in Theorem \ref{syzisosurf} as a corollary of
Theorem \ref{syzisocurv} and Proposition \ref{isofsyz}. Let us
also note that Theorem \ref{syzisocurv} and Theorem
\ref{syzisosurf} can easily be generalized to higher dimension and
the proofs are completely analogous to the ones given here. For
example, the $\mu$-basis of a curve in $\mathbb{P}^{m}$ consists
of $m-1$ syzygies whose degrees in $s$ and $\bar{s}$ sum up to
$d$. We are now ready to show our main result.

\newtheorem{impleqsurf}[hst]{Theorem}
\begin{impleqsurf}
Let $(q_1,q_2)$ be a $\mu$-basis of the parametrization
$\Phi_\mathcal{S} : \mathbb{P}^{1} \times \mathbb{P}^{1}
\dashrightarrow \mathbb{P}^{3}$. Then
$$Res(q_1,q_2)=F_\mathcal{S}^{deg(\Phi_\mathcal{S})}$$ where
$F_\mathcal{S}$ is an implicit equation of the ruled surface
$\mathcal{S}$ and where the resultant is taken with respect to $s$
and $\bar{s}$.
\end{impleqsurf}

\noindent \textit{Proof:} First, we can ensure that the hypotheses
of Proposition \ref{isofsyz} are fulfilled by performing a generic
linear coordinate change in $\mathbb{P}^{1} \times
\mathbb{P}^{1}$, which leaves both the implicit equation and the
resultant unchanged (up to multiplication by a constant). We will show that
$Res(q_1,q_2)$ is the power of an irreducible polynomial, i.e.
that it defines an irreducible hypersurface in $\mathbb{P}^3$. Let
us consider the incidence variety $\mathcal{W} := \{ ( \:
(s_0:\bar{s_0}),(x_0:y_0:z_0:w_0)) \in \mathbb{P}^1 \times
\mathbb{P}^{3} \: | \: q_i(s_0,\bar{s_0},x_0,y_0,z_0,w_0)=0 \: \}
$ then we have the following diagram

$$\begin{xy}
 \xymatrix{
    \mathcal{W}\ar@{->}[rrr]^{\pi_2} \ar@{->}[dd]_{\pi_1} & & & \mathbb{P}^3  \\
                                            & & &               \\
\mathbb{P}^{1}          & & &
    }
 \end{xy}$$ \\
where $\pi_1$ and $\pi_2$ are the canonical projections.
$\mathcal{W}$ is a vector bundle over $\mathbb{P}^{1}$, as the
$q_i$ are linear in $x,y,z,$ and $w$, and for any parameter
$(s_0:\bar{s_0})$ the fiber is a $\mathbb{K}$-vector space of
codimension 2 (because $q_1(s_0,\bar{s_0})$ and
$q_2(s_0,\bar{s_0})$ are linearly independent, as was proved in
\cite[Sect. 2, Prop. 3]{MR2015108}). As $\mathbb{P}^{1}$ is
irreducible, it follows that $\mathcal{W}$ is irreducible too (see
\cite[Ch.6, Th.8]{MR0447223}), hence so is $Im(\pi_2)$. (If
$Im(\pi_2)=A \cup B$ for two closed sets $A$ and $B$,
$\mathcal{W}=\pi_2^{-1}(A) \cup \pi_2^{-1}(B)$, which implies
$\mathcal{W}=\pi_2^{-1}(A)$ or $\mathcal{W}=\pi_2^{-1}(B)$, since
$\mathcal{W}$ is irreducible and, consequently, $Im(\pi_2)=A$ or
$Im(\pi_2)=B$). Now the points of $Im(\pi_2)$ are exactly those
for which the $q_i$ have a common zero in $s$ and $\bar{s}$, so by
definition of the resultant they are the zeros of $Res(q_1,q_2)$.
In other words, we have shown that $V(Res(q_1,q_2))=Im(\pi_2)$ is
irreducible, so $Res(q_1,q_2)$ is the power
of an irreducible polynomial.

 By definition, the syzygies of $\Phi_\mathcal{S}$ vanish
on all of $Im(\Phi_\mathcal{S})$ and hence on all of
$\mathcal{S}$, so $F_\mathcal{S} \: | \: Res(q_1,q_2)$. This
implies that $Res(q_1,q_2)$ is a power of $F_\mathcal{S}$ and it
remains to verify that it has the correct degree
$deg(\Phi_\mathcal{S})
\cdot deg(\mathcal{S})$.

 In the proof of Theorem \ref{isofsyz}, we have seen the
isomorphism of $R$-modules $$\begin{array}{rcl}
    \varphi : \qquad \textsl{Syz}_{R}(f_{0},f_{1},f_{2},f_{3}) & \rightarrow&  \textsl{Syz}(p_{03},p_{13},p_{23})
      \nonumber \\
     h_0x+h_1y+h_2z+h_3w &\mapsto& h_0x+h_1y+h_2z   \nonumber
     \end{array}$$
between the syzygies of the parametrization $\Phi_\mathcal{S}$ and
of the parametrization $\Phi_\mathcal{C}$ of its associated curve
$\mathcal{C}$. By abuse of notation, we will not differentiate between $\varphi$ and its extension to the morphism of
$R$-algebras $\varphi :  R[x,y,z,w] \rightarrow R[x,y,z]$ defined by $\varphi(x)=x$, $\varphi(y)=y$, 
$\varphi(z)=z$, and $\varphi(w)=0$.

As remarked earlier on, $\varphi$ takes $\mu$-bases to
$\mu$-bases, so $(\varphi(q_1),\varphi(q_2))$ is a $\mu$-basis of
$\Phi_\mathcal{C}$. Applying Theorem \ref{implcurv} yields
$$\begin{array}{rcl}
 F_\mathcal{C}^{deg(\Phi_\mathcal{C})} &=& Res(\varphi(q_1),\varphi(q_2))   \nonumber \\
         &=& \varphi(Res(q_1,q_2))  \nonumber
 \end{array}$$
where the last equality is true, because $\varphi$ is the
specialisation $w=0$ and as such commutes with the resultant.
Finally, we have $deg(\varphi(Res(q_1,q_2)))=deg(Res(q_1,q_2))$,
as $Res(q_1,q_2)$ is homogeneous, which shows that
$$deg(Res(q_1,q_2))=deg(\Phi_\mathcal{C}) \cdot deg(\mathcal{C})=deg(\Phi_\mathcal{S}) \cdot
deg(\mathcal{S})$$ so $Res(q_1,q_2)$ has indeed the correct
degree, which concludes the proof.
 \nolinebreak $\square$

\section{Algorithm and example}\label{algoex}

\noindent In this section, we give a detailed description of a new algorithm to 
compute a $\mu$-basis of a rational ruled surface based on the one-to-one
correspondence between the syzygies of a ruled surface and its associated curve: 
As we have remarked, a $\mu$-basis of the ruled surface can be obtained by computing
a $\mu$-basis of its associated curve (e.g. with the algorithm presented in \cite{ChenWangCurves})
and taking its image under the isomorphism \eqref{syzisom} in the proof of Proposition \ref{isofsyz}.
In particular, this method has the same computational complexity as the curve algorithm that
is used (since all the other steps in the algorithm are immediate), which makes it very efficient.

While it is convenient to work in the homogeneous setting for theoretical considerations,
actual computations should be done after dehomogenizing, i.e. setting $\bar{s}=1$ and $\bar{t}=1$ in
the parametrization \eqref{eqsurfdef}. In other words, we switch to the affine parametrization
$$\begin{array}{rcl}
    \Phi^{\texttt{aff}}_\mathcal{S} : \quad \qquad \mathbb{K}^{2}& \dashrightarrow&  \mathbb{K}^{3}
      \nonumber \\
     (s,t) &\mapsto& (\frac{f_0(s,t)}{f_3(s,t)},\frac{f_1(s,t)}{f_3(s,t)},\frac{f_2(s,t)}{f_3(s,t)})
\end{array}$$
where $f_i = f_{i0}(s)+t f_{i1}(s) \in \mathbb{K}[s,t]$. We remark that bihomogeneous polynomials of a fixed degree are in one-to-one correspondence to their dehomogenized counterparts and that this correspondence commutes with syzygy computations, resultants, etc. As a consequence, all the results in this paper are equally valid in the affine setting, so the $\mu$-basis and the implicit equation
can be obtained by computing their affine analogues and then re\-ho\-mo\-ge\-ni\-zing them.\\

\begin{quotation}
\noindent ALGORITHM ($\mu$-basis of a ruled surface)\\

\noindent INPUT: $f_i \in \mathbb{K}[s,t]$ for $i=0,1,2,3$ \\

\begin{enumerate}
 \item Check whether $deg_t(f_{i})=1$ for all $i$. If yes, set
$f_{i0}(s)=f_i(s,0)$ and $f_{i1}(s)=\frac{d}{dt}f_i(s,t)$ for all $i$. If not, return an error message.\\
 
 \item  Check whether $max(deg_s(f_{i1})) \geq max(deg_s(f_{i0}))$. If not,
        interchange $f_{i1}$ and $f_{i0}$ for all $i=0,1,2,3$.\\

 \item {Check whether $gcd(f_{30},f_{31})=1$. If not, check if there is $i \in \{0,1,2\}$ 
       such that $gcd(f_{i0},f_{i1})=1$.
       \begin{itemize}
 \item If there is such an $i$, interchange $f_i$ and $f_3$.

 \item If not, replace $f_3$ by $\alpha f_0 + \beta f_1 + \gamma f_2 + f_3$ for generic 
       $\alpha,\beta,\gamma \in \mathbb{K}$.\\
\end{itemize}}

\item Set $p_{i3}=f_{i0}f_{31}-f_{i1}f_{30}$ for $i=0,1,2$.\\

\item Calculate a $\mu$-basis
$(\tilde{q}_1,\tilde{q}_2)=(q_{11}x+q_{12}y+q_{13}z,q_{21}x+q_{22}y+q_{23}z)$
of the curve defined by $p_{03}$, $p_{13}$, and $p_{23}$ with an algorithm for planar curves.\\

\item Set $q_j=q_{j1}x+q_{j2}y+q_{j3}z-\frac{q_{j1}f_{00}+q_{j2}f_{10}+q_{j3}f_{20}}{f_{30}}$ for $j=1,2$.\\

\end{enumerate}
\noindent OUTPUT: A $\mu$-basis $(q_1,q_2)$ of the parametrization $\Phi^{\texttt{aff}}_\mathcal{S}$\\
\end{quotation}

\noindent Note that the second step of the algorithm may lead to a denser polynomial $f_3$ if a 
coordinate change is necessary, because the support of $f_3$ after such a change becomes the union of
the supports of the $f_i$. However, $f_0,f_1$ and $f_2$ are not changed, as we only have
to ensure the (relatively weak) condition $gcd(f_{30},f_{31})=1$ and do not need ``full'' genericity. 

Throughout the paper, we have considered a ruled surface to be given by a parametrization which has degree one in $t$. However, such a surface can also be defined by a parametrization of higher degree in $t$, so it would be interesting to give a criterion for when a given parametrization corresponds to a ruled surface and in this case to be able to replace it by another one which is linear in $t$.

\subsection*{Illustrative example} \noindent Let us consider the ruled surface $\mathcal{S}$ defined
by the polynomials $\tilde{f}_0 = s^2+t(s^2-1)$, $\tilde{f}_1=1+t(-s^2+1)$, $\tilde{f}_2= 1+t(-s^6+1)$, and 
$\tilde{f}_3=t(-s^6-2s^2)$. As $\tilde{f}_{30}=0$ and $\tilde{f}_3=s^2$ are not coprime, we interchange $\tilde{f}_3$ and $\tilde{f}_0$ and consider the new parametrization of $\mathcal{S}$
$$\begin{array}{rcl}
    f_0 & = & t(-s^6-2s^2)   \nonumber \\
    f_1 & = & 1+t(-s^2+1)   \nonumber \\
    f_2 & = & 1+t(-s^6+1)   \nonumber \\
    f_3 & = & s^2+t(s^2-1)   \nonumber
\end{array}$$ where $gcd(f_{30},f_{31})=1$.

\noindent Then its associated curve $\mathcal{C}$ is parametrized
by the Pl\"{u}cker coordinates
$$    p_{03}  =  s^8+2s^4    \qquad \qquad p_{13} =  s^4-1  \qquad \qquad    p_{23}  =  s^8-1$$

\noindent and we have $deg(\Phi_\mathcal{C}) \cdot
deg(\mathcal{C})=deg(\Phi_\mathcal{S}) \cdot deg(\mathcal{S})=8$
which follows from the degree formulae. Next we compute the following
$\mu$-basis for $\Phi_\mathcal{C}$ with a suitable algorithm:
$$\begin{array}{rcl}
    \tilde{q}_1 & = & (s^4+1)y-z   \nonumber \\
    \tilde{q}_2 & = & (-s^4+1)x-y+(s^4+1)z   \nonumber
\end{array}$$

\noindent Applying the isomorphism $\varphi^{-1}$ yields the following $\mu$-basis for
$\Phi_\mathcal{S}$ $$\begin{array}{rcl}
    q_1 & = & (s^4+1)y-z-s^2   \nonumber \\
    q_2 & = & (-s^4+1)x-y+(s^4+1)z-s^2   \nonumber
\end{array}$$

\noindent and we obtain
\begin{eqnarray}
   Res(q_1,q_2)  &=&  (4x^2y^2-4xy^3+y^4-4x^2yz+2xy^2z+x^2z^2+4xyz^2 \nonumber \\
                 & &  -2y^2z^2-2xz^3+z^4-x^2+xy+2y^2-xz-4yz+2z^2)^2 \nonumber
\end{eqnarray} which is the square of an implicit equation
$F_\mathcal{S}$ of $\mathcal{S}$. \\

\noindent  We have seen and used the equality
$deg(\Phi_\mathcal{C}) \cdot
deg(\mathcal{C})=deg(\Phi_\mathcal{S}) \cdot deg(\mathcal{S})$
between the surface $\mathcal{S}$ and its associated curve
$\mathcal{C}$. It is natural to ask whether
$deg(\mathcal{C})=deg(\mathcal{S})$ also holds. However, this is
not true in our example: we have $deg(\mathcal{C})=2$, but
$deg(\mathcal{S}) =4$. According, to the corollary to Proposition
\ref{isofsyz}, we would have had to perform a generic coordinate
change in order to ensure the equality of the degrees. \\

\noindent  Let us compare the $\mu$-basis method to some others.
In our example, $F_\mathcal{S}^2$ is obtained as a determinant of
a $8\times8$-matrix, the Sylvester matrix of $q_1$ and $q_2$.
After dehomogenizing our surface and homogenizing back to
$\mathbb{P}^2$ we can use approximation complexes to implicitize,
as in \cite{MR2172855}, and we obtain $F_\mathcal{S}^2$ as the
quotient of a $28 \times 28$-determinant by a $12\times
12$-determinant and an additional term that arises because we add
a non-complete-intersection base point when passing from $\mathbb{P}^1 \times \mathbb{P}^1$
to $\mathbb{P}^2$, which is by far not as efficient.
\\ \noindent Another possibility is to use the classical formula
$F_\mathcal{S}^2(w=1)=Res(f_0-xf_3,f_1-yf_3,f_2-zf_3)$ combined
with an efficient method to calculate the resultant such as
\cite{MR2004036}. $F_\mathcal{S}^2$ is obtained as the determinant
of $10 \times 10$-matrix, which is larger than our Sylvester
matrix and whose entries are themselves determinants of smaller
matrices.

\section{Remark on the reparametrization problem for ruled surfaces}

\noindent In the proof of Theorem \ref{implcurv} about the
implicit equation of a planar curve, we reduced the general case
to the proper case by reparametrizing the curve. If the field $\mathbb{K}$ is of
characteristic zero, we know by the theorem of Castelnuovo that there exists 
a proper reparametrization for any rational surface, i.e. there exists a commutative diagram
$$\begin{xy}
 \xymatrix{
    \mathbb{P}^1 \times \mathbb{P}^1 \ar@{-->}[rrr]^{\Phi_\mathcal{S}} \ar@{-->}[dd]_\psi & & & \mathbb{P}^3  \\
                                            & & &               \\
    \mathbb{P}^1 \times \mathbb{P}^1 \ar@{-->}[rrruu]_{\Phi'_\mathcal{S}}         & & &
    }
 \end{xy}$$ \\
where $\psi=(\sigma,\tau)$ is of degree $deg(\mathcal{S})$ and
$\Phi'_\mathcal{S}$ is a proper reparametrization of
$\mathcal{S}$. As far as we know, this problem is yet to be solved
algorithmically. However, \cite{MR2218599} gives a criterion for
the existence of a reparametrization of a rational surface such
that $\sigma=\sigma(s,\bar{s})$ depends only on $s$ and $\bar{s}$
and $\tau=\tau(t,\bar{t})$ depends only on $t$ and $\bar{t}$ and
proposes an algorithm for its computation if it exists. If we
restrict our attention to ruled surfaces we can also treat the
case where $\tau=( \bar{t} \alpha+ t \beta,\bar{t} \gamma+ t
\delta)$ with $\alpha,\beta,\gamma,\delta \in
\mathbb{K}[s,\bar{s}]$ such that $\alpha \delta - \beta \gamma
\neq 0$ by using the associated curve. So let us suppose that
there exists a reparametrization such that we can write
\begin{equation} f_i=\bar{t} (\alpha f'_{i0}(\sigma)+\gamma f'_{i1}(\sigma))+
t (\beta f'_{i0}(\sigma) + \delta f'_{i1}(\sigma))
\label{repasurf} \end{equation} for $i=0,\dots,3$, where the
$f'_{ij}$ define a proper parametrization $\Phi'_\mathcal{S}$ of
$\mathcal{S}$. We can deduce that
$deg(\psi)=deg(\sigma)=deg(\Phi_\mathcal{S})$, because $\tau$ is a
homography with respect to $t$. We have the
following identity $$p_i= \left| \begin{array}{cc} f_{i0} & f_{i1} \\
f_{30} & f_{31}
\end{array} \right| =\left| \begin{array}{cc} \alpha f'_{i0}(\sigma)+\gamma
f'_{i1}(\sigma) & \beta f'_{i0}(\sigma) + \delta
f'_{i1}(\sigma) \\
\alpha f'_{30}(\sigma)+\gamma f'_{31}(\sigma) & \beta
f'_{30}(\sigma) + \delta f'_{31}(\sigma)
\end{array} \right| = (\alpha \delta - \beta \gamma)
p'_i(\sigma)$$ from which we conclude that $\sigma$ yields a
proper reparametrization of the associated curve in the generic
case $deg(\Phi_\mathcal{S})=deg(\Phi_\mathcal{C})$. On the other
hand, any $\lambda(s,\bar{s})$ defining a proper reparametrization
$p_i=p''_i(\lambda)$ of $\mathcal{C}$ differs from $\sigma$ only
by a homography, so we can assume $\lambda = \sigma$, which
provides us with a (naive) method for calculating the
reparametrization: We compute $\sigma$ with a reparametrization
algorithm for curves such as in \cite{MR2218599} and consider
\eqref{repasurf} as a linear system of equations by comparing the
coefficients of the left hand side and the right hand side, where
we leave the coefficients of $\alpha, \beta, \gamma, \delta$ and
the $f'_{ij}$ undetermined. Then any solution of this system
defines a proper reparametrization of the ruled surface. However,
the systems are generally too large and further research is needed
to develop an efficient algorithmic solution to the
reparametrization problem.

\section{Acknowledgements}
\noindent The author was partially supported by the French ANR ``Gecko'' and the European 
project ACS nr. IST FET open 6413.

\end{document}